\documentclass[12pt]{article}
\usepackage{amsbsy,amsfonts,amsmath,amssymb,amsthm,graphicx}

\newtheorem{Theorem}{Theorem}[section]

\newtheorem{Lemma}[Theorem]{Lemma}

\newtheorem{Corollary}[Theorem]{Corollary}

\numberwithin{equation}{section}

\renewcommand{\emptyset}{\varnothing}

\newcommand{\0}{{\bf o}}

\newcommand{\cc}{{\rm c}}

\newcommand{\cH}{{\cal H}}

\newcommand{\bean}{\begin{eqnarray*}}
\newcommand{\eean}{\end{eqnarray*}}
\newcommand{\bea}{\begin{eqnarray}}
\newcommand{\eea}{\end{eqnarray}}

\newcommand{\toL}{\stackrel{L^1}{\longrightarrow}}

 \newcommand{\sgn}{\operatorname{sgn}}

\def\E{{\mathbb{E}}}
\def\Pr{{\mathbb{P}}}

\def\R{\mathbb{R}}

\def\Z{\mathbb{Z}}

\def\1{{\bf 1}}
\def\N{\mathbb{N}}
	\newcommand{\Po}{{\cal P}}
	\newcommand{\cQ}{{\cal Q}}
	
	\newcommand{\tQ}{\tilde{Q}}
	
	\newcommand{\cS}{{\cal S}}

	\newcommand{\tii}{\tilde{i}}
	\newcommand{\eps}{{\varepsilon}}

	\newcommand{\X}{{\cal X}}
	\newcommand{\cX}{{\cal X}}
	\newcommand{\cY}{{\cal Y}}

	\newcommand{\cC}{{\cal C}}
	
	\newcommand{\cV}{{\cal V}}

 \newcommand{\la}{\lambda}

 {\nopagebreak {\hspace*{\fill}\rule{7pt}{7pt}}\\ }

\title{Supercritical phase of the random connection model
}

\author{Mathew D. Penrose$^1$
\\
University of Bath
}

 \footnotetext{ Department of
Mathematical Sciences, University of Bath, Bath BA2 7AY, \\
United
Kingdom: {\texttt m.d.penrose@bath.ac.uk} }

\begin{document}
\maketitle

\begin{abstract}
	Given $d \in \N, \lambda >0$, the random connection model (RCM) in a region $A \subseteq \R^d$ is a graph with vertex set given by  a  homogeneous Poisson point process of intensity $\lambda $ in $A$, with  an edge placed between each pair $x,y$ of vertices with probability $\phi(\|x-y\|)$, where $\phi: \R_+ \to [0,1]$ is a nonincreasing finite-range connection function.  We show that if $d \geq 3$ and $\lambda$ is strictly supercritical for the RCM	in $\R^d$, then the RCM remains supercritical if it is restricted to a region $A$  of the form $\R^2 \times [-K/2,K/2]^{d-2}$, provided $K$ is sufficiently large. This is a continuum analogue of a well-known result of Grimmett and Marstrand for lattice percolation.  We prove this by adapting Grimmett and Marstrand's original proof; Faggionato and Hartarsky have also proved this recently by other means.

\end{abstract}

\section{Introduction and statement of results}

An important result in high-dimensional  percolation theory says that for
Bernoulli  site (or bond) percolation in $\Z^d$, $d \geq 3$,  
if the density of open sites (or bonds) is supercritical, it
remains supercritical when restricted to an infinite slab of
sufficiently large thickness. This result was proved by Grimmett and 
Marstrand \cite{GM90}. It is important, for example, because it enables one
to prove exponential decay for the probability of large finite
clusters, and finite correlation lengths, in the whole of the
supercritical phase. 

As Grimmett and Marstrand say, their methods and conclusions
`are valid for a large category of processes including many bond
and site processes on crystalline lattices'.
In the present paper we derive an analogous result for a model of
{\em continuum} percolation, namely the {\em random connection model}
in which a graph on  Poisson points in $\R^d$
is  created by connecting any two points with probability dependent
on how far apart they are (described in more detail below).
While the random connection model was introduced some time ago \cite{P91,MR},
there has been renewed interest recently in this model and generalizations
thereof \cite{P2,HLM,GHMM22} because
of their importance in general network modelling and telecommunications
modelling.

Although some continuum analogues of
Grimmett and Marstrand's 
result have previously been developed (discussed later), 
the random connection
model seems to require a more sophisticated (or at least, different)
adaptation of their methods
than do previous results of this sort in the continuum.

\pagebreak


We now describe our model in more detail.
Let $d \in \N$ and let
$\phi: [0,\infty) \to [0,1]$ be a 
nonincreasing function, known as the {\em connection function}.
Also let $p \in [1,\infty]$
and let $\|\cdot\|$ be the $\ell^p$ norm on $\R^d$.
Given a locally finite point set $\X \subset \R^d$, 
let
$G(\X,\phi)$ denote the graph with vertex set $\X$, where
for each $\{x,y\} \subset \X$, the
edge $xy$ is included with probability 
$\phi (\|y-x\|)$,
independently of
the other pairs.

Given $\la >0$, 
let $\cH_\la$ denote a homogeneous Poisson point process
of intensity $\la$ in $\R^d$ (viewed as a random set of points
rather than a measure).
Given also $s >0$, set $\Lambda_s:=[-s,s]^d$
and let $\cH_{\la,s} := \cH_{\la} \cap \Lambda_s$.
We are interested in the graphs $G(\cH_\la,\phi)$ and 
$G(\cH_{\la,s}, \phi) $, which are
known as the  {\em random connection model} \cite{MR}
and {\em soft random geometric graph} \cite{P2}
respectively, with connection function $\phi$.

Let $\cH^{\0}_{\lambda}$ denote the point process $\cH_{\lambda} \cup \{\0\}$,
where $\0$ is the origin in $\R^d$. Given Borel $A \subseteq \R^d$
with $\0 \in A$,
for $k \in \N$, let $\pi_k(\phi,\lambda,A)$ denote the probability that
the component of $G(\cH^{\0}_{\lambda} \cap A,\phi)$   containing
the origin is of
order $k$ (i.e., has $k$ vertices).
 The {\em percolation probability} in $A$, denoted
\index{percolation probability} $\theta(\phi,\lambda,A)$ is
 the probability that $\0$ lies in an
infinite component of the graph $G(\cH^{\0}_{\lambda} \cap A,\phi)$, that
is, 
$$
\theta (\phi,\lambda,A) : = 1 - \sum_{k=1}^ \infty \pi_k(\phi,\lambda,A).
$$
A standard coupling argument shows that $\theta(\phi,\la,A)$ is nondecreasing in $\la$.
The {\em critical value}\index{critical value}
 (continuum percolation threshold)
\index{continuum percolation threshold}
  $\lambda_\cc(\phi,A)$  
  is defined by
\bea
\lambda_\cc(\phi,A) := \inf \{\lambda>0: \theta(\phi,\lambda,A)>0\}.  
\label{lambdacdef}
\eea
It is known (see \cite{MR}) that 
	$0 < \lambda_\cc (\phi,\R^d) < \infty$,
	provided $0 < \int_{\R^d} \phi(\|x\|)dx
	< \infty$.

	Given $M \in (0,\infty)$,
	let $S_M := \R^2 \times [-M/2,M/2]^{d-2}$ (an infinite
	slab of thickness $M$).  We now state
	our main result, which is
	an analogue for the random connection model
	of the classic result of Grimmett and Marstrand \cite{GM90}
	for lattice percolation.
	\begin{Theorem}
		\label{t:GMphi}
		Suppose $d \geq 3$ 
		and $\sup\{r >0: \phi(r) >0\} \in (0,\infty)$.
		Then
	\begin{align}
		\lambda_c(\phi,\R^d) = \lim_{M \to \infty} \lambda_c(\phi,S_M).
		\label{e:GMphi}
	\end{align}
	\end{Theorem}
	This result has also been proved recently by
	different means in independent (but prior) work of
	Faggionato and Hartarsky \cite{FH25}, which we shall discuss 
	further at the end of this section.

	As discussed in \cite{GM90},
	Grimmett and Marstrand's result has many applications in lattice
	percolation theory, and we would expect Theorem \ref{t:GMphi}
	to have a similar range of applications for the random connection
	model. In Section \ref{s:Giant} 
	we discuss just one such application, to the {\em giant
	component} phenomenon in soft random geometric graphs.


Another popular continuum percolation model is the (spherical)
{\em Boolean model},
where balls of fixed radius or i.i.d. random radii are centred on the points 
of $\cH_\la$, and one investigates the components of the resulting union of balls. For the Boolean model, an analogous result to \eqref{e:GMphi}
was given by Tanemura \cite{Tanemura} for fixed radius and
extended to bounded random radii
by Last, Penrose and Zuyev
\cite[Theorem 3.7]{LPZ} (see also \cite{Fag21} for a further generalization). 
Since balls of fixed radius $r$ overlap if and only if their centres
are distant at most $2r$ from each other, percolation in the Boolean model with
balls of fixed radius $r$ is equivalent to percolation in the random 
connection model with $\phi = {\bf 1}_{[0,2r]}$.
Therefore  the special case 
$\phi = {\bf 1}_{[0,2r]}$ of our result
\eqref{e:GMphi} was already known
due to Tanemura's result.


In all of \cite{Tanemura,LPZ,Fag21} the continuum model
is dealt with by a fine discretization of $\R^d$ into small
cubes, allowing the authors to directly apply a version of Grimmett
and Marstrand's result for site percolation, albeit 
on a more complicated
lattice than the usual nearest-neighbour lattice $\Z^d$. In the present
instance of the random connection model with more general $\phi$,
it is not clear that such a discretization argument works, since
the cubes in the discretization can be multiply occupied, meaning
that no convenient site-bond percolation approximation 
seems to be available.
Our proof  instead uses a coarser discretization
into unit cubes for adapting some of the counting arguments
in \cite{GM90}, while at the same time staying in the continuum for 
the proof as a whole.

Our method should work for a slightly wider class of connection
functions $\phi$. Provided there exists $\delta >0$
such that $\phi(r) > \delta$ for $0 < r < \delta$,
by only minor changes to the proof
we should be able to drop the assumption that
$\phi$ is nonincreasing. However the assumption that
$\phi$ has finite range is essential to our proof, so
finding a way to drop that assumption would be 
of great interest. Another
possible generalization would be to allow for $\cH_\la$
to be a {\em marked} Poisson process, with the connection
function allowed to depend on the marks of the endpoints.
This generalization has been considered in some of the recent
literature on this kind of model
(see e.g. \cite{Fag21,GHMM22,K25})
but we do not address it here.

Shortly before the first version of the present paper was completed,
Faggionato and Hartarsky proved a result that essentially implies
our Theorem \ref{t:GMphi} \cite[eq. (37)]{FH25}.
In fact, their work allows for $\cH_\lambda$ to be a marked point
process, as well as for a somewhat more general class of connection
functions (though still restricted to finite-range connections).
Their method relies on the 
`seedless' renormalization scheme  of Duminil-Copin, 
Kozma and Tassion \cite{DCKT}; the present paper shows that
at least in the unmarked case, it is also possible to adapt
the original argument of Grimmett and Marstrand.

\section{Preliminaries and lemmata}
\label{chgiant1}
\label{seccoper}

We introduce some further notation that will be used in proving
Theorem \ref{t:GMphi}.

For $A \subseteq \R^d$ and $x \in \R^d$, set $x+ A
:= \{a+x:a \in A\}$. Also  set
$\Lambda_r(x) := x + \Lambda_r$.

Given $A,B \subset \R^d$, and locally finite $\X \subset \R^d$,
we write $\{ A \leftrightarrow B $ in $G(\X,\phi)\}$
for the event that there exist $x \in \X \cap A $ and $y \in \X \cap B $
such that there is a path in $G(\X,\phi)$ from $x $ to $y$.
We write $\{ A \leftrightarrow \infty $ in $G(\X,\phi)\}$
for the intersection over all $n \geq 1$ of events
 $\{ A \leftrightarrow \R^d \setminus \Lambda_n $ in $G(\X,\phi)\}$.

 Given nonincreasing $\phi:[0,\infty) \to [0,1]$, set
 $\rho_\phi := \sup\{r  > 0:\phi(r) >0\}$, the {\em range} of $\phi$. 
 We show next that it suffices to prove Theorem \ref{t:GMphi}
 in the case where $\rho_\phi =1$, by a 
 {\em scaling argument}.
 Given $a >0$, set $\phi_a (r) :=  \phi(ar)$
 for all $r \geq 0$. 
 We first claim that for any $\lambda, a >0$ and
 Borel $A \subseteq \R^d$ with $\0 \in A$, setting
 $a^{-1}A:= \{a^{-1}x: x \in A\}$   we have
 \begin{align}
	 \label{e:scaling}
	 \theta (\phi,\lambda,A) = \theta (\phi_a, a^d \lambda, a^{-1}A).
 \end{align}
 Indeed, by the Mapping theorem \cite[Theorem 5.1]{LP}, 
 the image of $\cH_\la \cap A$
 under the mapping $x \mapsto a^{-1}x$ has the same distribution
 as $\cH_{a^d \la} \cap a^{-1} A$, and \eqref{e:scaling}
 follows.
 By \eqref{e:scaling} and 
 \eqref{lambdacdef},  
 $$
 \lambda_c(\phi,A) = \inf \{\lambda:\theta(\phi_a,a^d \lambda,a^{-1}A) 
 > 0 \} = a^{-d} \lambda_c(\phi_a, a^{-1} A).
 $$
 Therefore, if the conclusion \eqref{e:GMphi} of Theorem \ref{t:GMphi}
 holds for $\phi$, then also
 \begin{align*}
	 \lambda_c (\phi_a,\R^d) & = a^d \lambda_c (\phi,\R^d)
	 \\
	 & = 
	 a^d \lim_{M \to \infty} \lambda_c(\phi,S_M)
	 \\
	 & = \lim_{M \to \infty} \lambda_c ( \phi_a, a^{-1} S_M) 
	  = \lim_{M \to \infty} \lambda_c ( \phi_a,  S_{M/a}), 
 \end{align*}
 so \eqref{e:GMphi} holds for $\phi_a$ as well as for $\phi$.
 Since for any $\phi $ with $\rho_\phi \in (0,\infty)$,
 setting $a = \rho_\phi$ and $\phi' = \phi_{a^{-1}}$, we have
 $\phi = \phi'_a$,
 and $\rho_{\phi'} =1$, it suffices to prove 
 Theorem \ref{t:GMphi} in the case with $\rho_\phi =1$.
 Therefore,
 from now on until the end of Section \ref{s:continuing},
 we assume that $d \geq 3$ and  
$\phi: [0,\infty) \to [0,1]$ is nonincreasing
with $\rho_\phi :=
 \sup \{r>0: \phi(r) > 0\} =1$. 

We shall view $\phi$ as fixed from now on, and
for any locally finite $\X \subset \R^d$ 
write simply $G(\X)$ instead of $G(\X,\phi)$.
For $x \in \X$,
let $\cC_x(\X)$ denote the vertex set of the
component of $G(\X)$ containing $x$ (these components are 
often called {\em clusters}).  
Also given finite  $\cY \subset \R^d$, 
we write $\cC_\cY(\X)$ for $\cup_{y \in \cY } \cC_y(\X \cup \cY)$
(here $\cY$ might or might not be a subset of $\cX$).
We shall refer to $\cC_\cY(\X)$ as a
{\em set-cluster};  it is a finite union of clusters, associated
with the point set $\cY$.

For $i \in \{1,2,\ldots,d\}$,
let $e_i$ denote the $i$th coordinate unit vector, e.g.
$e_1 := (1,0,\dots,0)$.


We say a real-valued function $f$, defined on graphs
$(\cV,E)$,  with
$\cV \subset \R^d$ locally finite, is {\em increasing}  if $f(\cV,E) \leq
f(\cV',E')$
whenever $\cV \subset \cV'$ and $E \subset E'$. We say $f$ is {\em decreasing}
if $-f$ is increasing. Given $\la >0$ and given $\phi$,
we say an event  ${\cal E}$ 
on $G(\cH_\la)$
is {\em increasing} (resp. decreasing) 
if ${\bf 1}_{\cal E} $
is an increasing (resp. decreasing)
function of $G(\cH_\la)$.
\begin{Lemma} [Harris-FKG inequality]
	\label{lemFKG}
	Suppose $f,g$ are   measurable bounded  increasing real-valued
	functions
	defined on graphs $(\cV,E)$ with $\cV \subset \R^d$ locally finite.
	Then $$\E[f(G(\cH_{\la})) g(G(\cH_{\la})) ] \geq \E[f(G(\cH_{\la}))]
	\E[g(G(\cH_{\la}))]. $$
	The same inequality holds if $f$ and $g$ are both decreasing.
\end{Lemma}
\begin{proof}
See \cite{HLM}, where measurability issues are also dealt with.
\end{proof}

	Given $m, n \in \N$,
following \cite{GM90} we define the sets 
\begin{align}
	T(n) &:= \{x 
\in ( \tfrac{1}{2} +\Z)^d:
x \cdot e_1 =n- \tfrac12, 0 < x \cdot e_j < n
	{\rm ~for~} j=2,3,\ldots,d \} ;
\nonumber
	\\
	T(n,m) & := [n,n+m] \times [0,n]^{d-1},
	\label{e:Tnm}
\end{align}
and similarly
\begin{align*}
	T_2(n) := \{
		x \in ( \tfrac{1}{2} +\Z)^d:
	x \cdot e_2 =n- \tfrac12, 
	0 <  x \cdot e_j <  n,
	{\rm ~for~} 
	j=1,3,\ldots,d
\} ;
\nonumber
	\\
	T_2(n,m)  := 
	\{	x \in \R^d : n \leq x \cdot e_2 \leq n +m, 
	0 \leq  x \cdot e_j \leq  n, {\rm ~for~} j = 1,3, \ldots, d
	 \} .
\end{align*}
Next we adapt \cite[Lemma 4]{GM90}. For finite $A $ we write
$\#A$ for the number of elements of $A$.
\begin{Lemma}
	\label{l:likeGML4}
	Let $\eps >0$, $\lambda > \lambda_c(\phi,\R^d)$.
	For $r, k \in \N$ with $k >r$, 
	define the random set
	\begin{align}
		V(r,k) := \{x \in T(k): \Lambda_r \leftrightarrow
		\Lambda_{1/2}(x) {\rm ~in~} G(\cH_\la \cap 
		\Lambda_{k} 
		)\}.
		\label{e:defV}
	\end{align}
	Then there exists $m \in  \N$ with $m \geq 9$ such that
	for all $\ell \in \N$, we have
	\begin{align}
		\limsup_{k \to \infty} \Pr[\#V(m,k) < \ell] <
		 \eps.
		 \label{e:GML4}
	\end{align}
\end{Lemma}
\begin{proof}
	Since we assume $\lambda > \lambda_c(\phi,\R^d)$,
	we have $\theta(\phi, \lambda,\R^d) >0$,
	so $\Pr[N_\infty>0] >0$,
	where $N_\infty$ denotes the number of unbounded
	components of $G(\cH_\la)$. Hence by
	the zero-one law for the random connection model
	\cite[Theorem 2.14 and Proposition 2.6]{MR},
	$\Pr[N_\infty > 0] =1$.

	Set $w := d2^d$.
	Defining $F_n 
		:= \{\Lambda_n \leftrightarrow \infty\
		{\rm ~in~} G(\cH_\la) \}^c $, 
		we have $\Pr[ \cap_{n=1}^\infty F_n] = \Pr[N_\infty =0]
		= 0$. Therefore we can and do
		choose $m \in \N$ such that $m \geq 9$ and
		$\Pr[F_m ] < \eps^w$.

	Given $k \in \N $ with $k > m$,
	set 
	$$
	U_k:= \{x \in (\tfrac12 + \Z )^d 
	\cap \Lambda_k \setminus \Lambda_{k-1}:
	\Lambda_m \leftrightarrow\Lambda_{1/2}(x) 
	{\rm ~ in~} G(\cH_\la \cap \Lambda_{k}  )\}.
	$$
	Then, as in the proof of \cite[Lemma 4]{GM90},
	there is a collection of $w$ decreasing
	events, denoted $A_1,\ldots,A_w$
	say, each having the same probability as the event
	$\{\# V(m,k) < \ell \}$, such that 
	$$
	\{ \# U_k < w \ell \} \supseteq \cap_{i=1}^w A_i.
	$$
	Hence by the Harris-FKG inequality (Lemma \ref{lemFKG})
	\begin{align}
	\Pr[\# U_k  < w \ell ] \geq \Pr[\# V(m,k) < \ell]^w.
		\label{e:forGML4}
	\end{align}

	Now let
	$	\Delta_k U_k  :=  
		 ( \cup_{z \in U_k} \Lambda_{3/2}(z))
		\setminus \Lambda_{k}. $
		By the union bound, $\Delta_k U_k$ has Lebesgue $d$-measure
		at most $3^d \# U_k$, so that
	\begin{align*} 
		\Pr[U_{k+1}= \emptyset | 1 \leq \# U_k < w \ell ]
		& \geq \Pr[ \cH_\lambda \cap \Delta_k U_k = \emptyset
		| 1 \leq \# U_k < w \ell ]
		\\
		& \geq \exp (- 3^d \lambda w \ell).  
	\end{align*} 
	Also if $\Lambda_m \leftrightarrow \infty$ in $G(\cH_\la)$ then
	$U_k \neq \emptyset$.
	Hence
	\begin{align*}
		e^{-3^d \lambda w \ell} \Pr[\# U_k < w\ell,
		\Lambda_m \leftrightarrow \infty {\rm ~ in ~} G(\cH_\la)] & 
		\leq
		e^{-3^d \lambda w \ell} \Pr[ 1 \leq \# U_k < w \ell]
		\\
		&	\leq \Pr[ 1 \leq \# U_k < w \ell, U_{k+1} = \emptyset]
		= \Pr[ E_k],
	\end{align*}
	say, where the events $E_k, k > m$
	are disjoint so their probabilities are summable in $k$.
	Therefore $\Pr[\# U_k < w\ell,
		\Lambda_m \leftrightarrow
		\infty {\rm ~ in~} G(\cH_\la)] \to 0$ as $k \to \infty$.
		Thus, since
	\begin{align*}
		\Pr[ \# U_k < w \ell ] \leq 
		\Pr[ \# U_k < w \ell, \Lambda_m \leftrightarrow \infty
		{\rm ~ in~} G(\cH_\la) 	] +
		\Pr[F_m],
	\end{align*}
	using our choice of $m$ at the start of this proof
	we have
		$$
		\limsup_{k \to \infty} \Pr[\#U_k < w \ell] \leq  \Pr[F_m]
		< \eps^w,
		$$
		and since  \eqref{e:forGML4} gives us
		$
		\Pr[\# V(m,k) < \ell ] \leq
		\Pr[\# U_k < w \ell]^{1/w}, 
		$
		we can deduce
		\eqref{e:GML4}.
\end{proof}

Next, we shall give an analogue to \cite[Lemma 2]{GM90}. 
Given $\eta \subset \R^d$, 
set
\begin{align}
\overline{\eta} : = \cup_{z \in (\tfrac12 + \Z  )^d: \Lambda_{1/2}(z) \cap \eta
\neq \emptyset} \Lambda_{1/2}(z).
	\label{e:bardef}
\end{align}
In other words, $\overline{\eta}$ is the union
of all cubes of the form $\Lambda_{1/2}(z),  
z \in ( \tfrac12 + \Z)^d$ having non-empty intersection with
$\eta$.  Given also $ B \subset \R^d$ with $B = \overline{B \cap (\frac12 \cap
\Z)^d}$,
define the sets
\begin{align*}
	\Delta(\eta;B) & := \{z \in (\tfrac12 + \Z)^d
	\cap B \setminus \overline{\eta} 
	: \Lambda_{1/2}(z) \cap \overline{\eta 
	\cap B} \neq \emptyset\}; 
\\
	\eta^{+}_B & := (\eta \cap B) \cup \Delta(\eta;B) .
\end{align*}
In other words, $\Delta(\eta;B)$ is the external boundary
(in the $\ell^\infty$ sense) of
$\eta$ with respect to $B$ in the discretization. For later
use we also define
\begin{align}
	\Delta^2(\eta;B)  := \Delta (\eta^+_B;B);
	~~~~~~
	\eta^{++}_B  := 
	\eta_B^+ \cup \Delta^2(\eta;B).
	\label{e:++def}
\end{align}

\begin{Lemma}
	\label{l:likeGML2}
	Let $\lambda >0$ and $R,  B \subset \R^d$ with
	$ B = \overline{B \cap (\frac12 + \Z)^d}$.
	Let $K$ be a random subset
	of $( \frac12 + \Z)^d \cap B$, independent of 
	$\cH_\la$, 
	with $R^{+}_B 
	\cap \overline{K} = \emptyset$ almost surely.
	Define the random set 
	$$
	U = \{ z \in \Delta(R;B): \overline{\{z\}}
	\leftrightarrow \overline{K} {\rm ~in~}
	G( \cH_\la \cap  B \setminus 
	\overline{R}) 
	\}.
	$$
	Then 
	for each $t \in \N$,
	\begin{align}
	\Pr[\# U \leq t ]  
	\leq e^{3^d \lambda t}
		\Pr[\{\overline{R}
		\leftrightarrow \overline{K}
	{\rm ~ in~} G( \cH_\la  \cap B )
		\}^c].
		\label{e:forGML2}
	\end{align}
\end{Lemma}
\begin{proof}
%
	Let $\alpha \subset \Delta(R;B)$. If $U= \alpha$
	and $\cH_\la \cap ( \cup_{z \in \alpha} \Lambda_{3/2}(z)) \cap
	\overline{R \cap B} = \emptyset$, then event
	$\{\overline{R} \leftrightarrow \overline{K} {\rm ~in~}
	G(\cH_\la \cap B)\}$ does not occur, since
	any path from $\overline{K}$ to $\overline{R}$
	must include a step from a vertex $v \in \overline{U}$
	to  a vertex $w \in \overline{ R \cap B} $, and
	if $v \in \overline{\{z\}}$ say for some $z \in \overline{U}$,
	then $w \in \Lambda_{3/2}(z)$ by the triangle inequality.
	By the independence properties
	of $\cH_\la$ and the independence assumption on $K$,
	the events $\{U= \alpha\}$  
	and $\{\cH_\la \cap ( \cup_{z \in \alpha} \Lambda_{3/2}(z)) \cap
	\overline{R \cap B} = \emptyset \}$ are independent, and 
	$\cup_{z \in \alpha} \Lambda_{3/2}(z)$ has
	Lebesgue $d$-measure at most $3^d\# \alpha$. Therefore
	\begin{align*}
		\Pr[ \{\overline{R} \leftrightarrow \overline{K} {\rm ~in~}
		G(\cH_\la \cap B)\}^c ] 
		& \geq \sum_{\alpha \subset \Delta(R;B)} \Pr[U= \alpha]
		\Pr[ \cH_{\la}\cap ( \cup_{z \in \alpha }
		\Lambda_{3/2}(z))  \cap \overline{R \cap B} = \emptyset ]
		\\
		& \geq \Pr[\# U \leq t] \exp(- 3^d \lambda t),
	\end{align*}
%
	and hence
	\eqref{e:forGML2}.
\end{proof}

\section{Algorithm: First step}
\label{s:step0}

	It is not hard to see that $\lambda_c(\phi,S_M)$ is nonincreasing in $M$
	and $\lambda_c(\phi,S_M) \geq \lambda_c(\phi,\R^d)$ for all $M$,
	so to prove  \eqref{e:GMphi}
	it suffices to prove that if $\mu > \lambda_c(\phi,\R^d)$
	then we can find $M>0$ such that $\mu \geq \lambda_c(\phi,S_M)$.
 We therefore fix arbitrary $\mu
 > \lambda > \lambda_c(\phi,\R^d)$.
We shall show for some $M >0$ (depending on $\lambda $ and $\mu$)
that $\0 \leftrightarrow
\infty {\rm ~in~} G(\cH_\mu^\0 \cap S_M)$ with positive probability,
i.e. $\theta(\phi, \lambda, S_M) >0$, so that
$\lambda_c(\phi,S_M) \leq \mu$
as required.

Our proof will be based 
on an algorithm which will 
generate a nondecreasing sequence  of subsets of $\cH_{\mu}^\0$,
denoted $\xi_1, \xi_2, \xi_3, \ldots$, such that
$\0 \in \xi_i$ and 
the subgraph of $G(\cH^{\0}_{\mu})$ induced by $\xi_i$
(denoted $G(\xi_i)$ for short)
is connected for each $i \in \N$.
Since $G(\xi_i)$ is connected and therefore part of 
 a single cluster of $G(\cH_{\mu}^{\0})$, 
we shall refer to $\xi_i$ as the {\em sub-cluster at the origin}
at the end of the $i$th step of the algorithm. 

The general idea of the algorithm goes as follows.
We shall think of $\cH_\la$ as being a subset of $\cH_\mu$, and
take 
$\xi_1$ to be the restriction of $\cH_\la$ to $\Lambda_m$, where
$m$ is a  large fixed integer to be chosen in a manner we shall describe later,
together with the point $\0$ and some of the points of $(\cH_\mu \setminus
\cH_\la) \cap \Lambda_1$ (we will spell out which ones later on).
We shall then consider 
a sequence of non-homogeneous Poisson point processes
$\cS_2$, $\cS_3$, $\cS_4, \ldots$,
each of which can be thought of as a sub-process
of $\cH_\mu$. Each point process $\cS_n$ will be
restricted to a bounded spatial region $A_n$, 
where the $A_n$ are generated sequentially
and for each $n$,
$A_n$ depends on the history of the algorithm prior to the
$n$th
step.
To get from $\xi_n$ to $\xi_{n+1}$ we
determine the set-cluster $\cC_{\xi_n}(\cS_{n+1})$,
and  add the points of this set-cluster.

Of use to us is the following {\em  sequential construction}
of $\cC_{\xi_n}(\cS_{n+1})$.
%
Denote the points of $\xi_n$ 
as {\em active points} and let $g(\cdot)$  be the initial  intensity function
of {\em unexplored  points}, i.e.
the intensity function of $\cS_{n+1}$, representing  Poisson points that
are not yet generated.

Next, choose an active point $x$ and generate a Poisson process
of intensity 
$ g(\cdot) \phi(\|\cdot -x\|) $,
representing the
previously unexplored points
of $\cS_{n+1}$ that are connected directly to $x$. Label all the new 
points as `active', and change the status
of $x$ from `active' to `finished'. Also change the intensity of
unexplored points of $\cS_{n+1}$ from $g(\cdot)$ to
$g(\cdot) (1- \phi(\|\cdot -x\|))$.

Then pick a new active point and repeat the above, using the new intensity of
unexplored points. Keep repeating until we run out of active points,
then stop.

We shall refer to the above procedure as
{\em growing the cluster sequentially}. This method is described in detail
for the case $\phi = {\bf 1}_{[0,1]}$ 
in \cite{P96}, and for the general random connection
model in \cite{MPS}.

The main point of this construction is that
the region  distant further than 1 (in our  norm $\|\cdot\|$) from 
$\cC_{\xi_n}(\cS_{n+1})$
remains unexplored at the end of this step, and we
can use the whole of our Poisson process $\cH_\mu$ in
this unexplored region, in later steps.

A slightly different (and possibly more elementary)
way to perform  the sequential 
construction, which  would also work for our purposes, goes as follows. Divide
$\R^d$ into unit cubes $B_i$ centred on the points of
$(\frac12 + \Z)^d$. At each stage, given $\xi_n$, 
the  next region $A_{n+1}$ will be a finite union of
some of these cubes. Let us call $\xi_n$
the current sub-cluster at the origin. As well as
$\xi_n$ there will be a second point process $\varphi_n$
of points that have been revealed but
are not connected to $\xi_n$.
Next, pick one of the cubes $B_i$ that is contained
in $A_{n+1}$ and within $\ell^\infty$ distance 1 of
at least one point of the current sub-cluster at the origin. 
Reveal the points of the Poisson process $\cS_{n+1} \cap B_i$
(new points), and the status (present/absent) of all potential
edges between the new points and between new points and old points.
This may cause some of the new points  (and also possibly
in some of the previously revealed points that were not previously
in the sub-cluster at the origin)
to be added to the current sub-cluster at the origin.
Then pick another unit cube $B_i$ that is contained
in $A_{n+1}$ and within unit $\ell^\infty$ distance
of the current sub-cluster at the origin. Repeat
until no such $B_i$ exists. Then conclude this step;
at the end of this step
none of the cubes $B_i \subset A_{n+1}$ that
are $\ell^\infty$ distance more than $1$ from 
$\xi_{n+1}$ will have been explored.

Let $\eps_1 := (\mu- \lambda)/40$.
By the Superposition theorem \cite[Theorem 3.3]{LP}, we can
(and do) think of $\cH_{\mu}$ as being decomposed
into a homogeneous Poisson process $\cH_\la$ of intensity
$\lambda$ in $\R^d$, along with 40 independent homogeneous Poisson
processes (which we think of as `sprinkled points')
of intensity $\eps_1$.
The algorithm will  explore these Poisson processes in
various regions of space that are determined by the prior
history of our exploration. 
That is, the point process $\cS_{n+1}$ will be the
restriction of $\cH_\la$ to
some bounded region $A^{(1)}_{n+1}$ together with the restriction
of two of the point processes of sprinkled points
to some other region $A_{n+1}^{(2)}$.
The non-homogeneity of $\cS_{n+1}$ arises from the fact
that 
$A_{n+1}^{(1)} \neq A_{n+1}^{(2)}$ 
in general.
We shall use the sprinkled points
to cross certain boundaries that arise between explored and
unexplored regions of $\R^d$.

We now define our version of what is called a {\em pad}
in \cite{GM90} or a {\em seed} in \cite{Grimmett}.
Given $x \in \Z^d$ and $m \in \N$, we divide the cube $\Lambda_m(x) $
into unit subcubes. 
We say that $\Lambda_m(x)$
is a {\em seed} if
each of these subcubes contains at least one point of $\cH_\la$,
and moreover
$G(\cH_\la \cap \Lambda_m(x))$ is connected.
Note that we require our seeds to be centred on points of $\Z^d$.
%

Let $\eps_0 = 1/9999$. We shall explain our choice of 
this parameter later on.
We  choose $m$, and introduce various other constants,
as follows. First we define
\begin{align}
	\delta_1 := e^{- 5^d \eps_1} (\eps_1/(2d)^d)^{(10d)^d}
	\phi((d+1)/(2d))^{(10d)^d-1 } 
	\phi(1/2)^2 .
	\label{e:delta1}
\end{align}
Next, define
\begin{align}
	t_1 : = 9^d \lceil (\log (\eps_0^{-1}))/\delta_1\rceil.
	\label{e:tdef}
\end{align}
By Lemma \ref{l:likeGML4}, we can and do choose $m \in \N$
with $m \geq 9$
such that for any $\ell \in \N$,
	\begin{align}
		\limsup_{k \to \infty} \Pr[\# V(m,k) \leq \ell ] \leq 
	 e^{-3^d \lambda t_1} \eps_0/3. 
		\label{e:V'0}
	\end{align}

	With $m$ (and $\lambda$) now fixed, we  further define
\begin{align}
	\eps_2 & := 
\Pr[\Lambda_m ~{\rm is~a~seed}] >0;
	\label{e:eps2}
	\\
	\eps_3 & := e^{-2\eps_1} (\eps_1/(2d)^d)^{2(2d)^d} 
	\phi((d+1)/(2d))^{2(2d)^d-1} 
	\phi(1/2)^{2}, 
	\label{e:defeps3}
\end{align}
and then define
	\begin{align}
	t_2 
		: = \lceil (2m)^{d-1} \eps_3^{-1}\eps_2^{-1} (3^d \lambda t_1
	 + \log (2/\eps_0)) \rceil.
		\label{e:ell}
	\end{align}
Finally, using \eqref{e:V'0} we can and do take
$n \in \N$ to be a multiple of $2m$ that is sufficiently large that
\begin{align}
	\Pr[\# V(m,n) < t_2] 
		\leq e^{-3^d \lambda t_1} \eps_0/2.
		\label{e:ndef}
\end{align}
	The constants $\eps_0, \eps_1, \eps_2, \eps_3, \delta_1,
	t_1, t_2, m$ and $n$
	will remain fixed for the rest of the proof.

	In the first step of the algorithm we sample
	the Poisson process $\cH_\la \cap \Lambda_m$, and
	also an independent homogeneous Poisson process
	of intensity $\eps_1$
	on $\Lambda_{1/2}$, which we denote $\Po_1$,
	and also the status (present/absent) of all potential
	edges between the points thus found.

	Given any finite set $\eta \subset \R^d$ 
	and any $z \in (\frac12 + \Z)^d$,
	let $\eta\|z$ denote the first point of
	$\eta \cap \Lambda_{1/2}(z)$ in
	the lexicographic ordering if $\eta \cap \Lambda_{1/2}(z) 
	\neq \emptyset$; otherwise set $\eta \| z := \0$.

	Let $E_1$ be the event that (i)
	$\Lambda_m $ is a seed, and
	(ii) $\0$ is connected through a path in 
	$G(
	\{\0, \cH_{\la}\|\0 \}
	\cup \Po_1 
	)
	$ 
	to the point $\cH_{\la}\|\0$.
	Then we claim  $\Pr[E_1 ] \geq \eps_2 \delta_2 >0 $, where we set
	$$
	\delta_2 := e^{-\eps_1} (\eps_1/(2d)^d)^{(2d)^d}
	\phi((d+1)/(2d))^{(2d)^d-1} \phi(1/2)^2.
	$$
	Indeed, $\delta_2$ is
	a lower bound for the probability that if we divide
	$\Lambda_{1/2}$ into subcubes of side $1/(2d)$,
	each subcube contains exactly one point of $\Po_1$,
	and the graph $G(\Po_1)$ is connected, and both the origin
	and $\cH_\la \| \0$
are connected
	to the point of $\Po_1 $ in their subcube.
Note that any two points in neighbouring subcubes
are distant at most $(d+1)/(2d)$ from each other in our
	chosen $\ell^p$ norm, while
	any two points in the same subcube are
	distant at most $\frac12$ from each other. We shall make
	a number of similar estimates later on.

If $E_1$ occurs, we say 
the first step of the algorithm is successful,
and set $\xi_1 : =  \Po_1 \cup
(\cH_\la^{\0} \cap \Lambda_m)$.  If $E_1$ does not
 occur set $\xi_n = \{\0\} $ for all $n \in \N$.

 \section{Algorithm: Second step}
\label{s:step1}

 We now describe the next step of the algorithm.
 Assume $E_1$ occurs. Then $\Lambda_m$ is a seed and
 $\overline{\xi_1} = \Lambda_m$, where we recall the notation
 $\overline{\eta}$ from \eqref{e:bardef}.

 Let $\Po_2$ and $\cQ_2$ be  homogeneous Poisson processes
 of intensity $\eps_1$ in $\Lambda_{n+1}$, independent
 of $\cH_\la,$ $\Po_1$ and each other.
Then, recalling the definition of $T(n,m)$, $T_2(n,m)$ at 
	\eqref{e:Tnm}, define the set 
\begin{align}
	\xi_2 & := 
	\cC_{\xi_1} (
	\Po_2 \cup \cQ_2 \cup (\cH_\la \cap [( \Lambda_{n} \setminus \Lambda_{m+1} )
	\cup	T(n,m) \cup T_2(n,m)] ) ).
	\label{e:defxi1}
\end{align}
Then $\xi_2$
is the new current value of the sub-cluster
at the origin, at the end of this step of the algorithm.

We denote by $E_2$ the event that there is a seed contained in 
$\overline{\xi_2} \cap T(n,m)$, and by $E'_2$ the event that
there is a seed contained in $\overline{\xi_2} \cap T_2(n,m)$.
If $E_2 \cap E'_2$ occurs, then we say that this step of the
algorithm  is successful.  In the next lemma we estimate
$\Pr[(E_2 \cap E'_2)^c | E_1]$.

Suppose this step is successful. Then
there is a seed contained in $T(n,m) \cap \overline{\xi_2}$.
Let $v_1$ be the centre of the first such seed in the lexicographic
ordering. 
There is also a seed contained in $T_2(n,m) \cap \overline{\xi_2}$.
Let $v_2$ be the centre of the first such seed in the lexicographic
ordering. We say the seeds $\Lambda_m(v_1)$ and
$\Lambda_m(v_2)$ have been {\em identified} at the end
of this step.

\begin{Lemma}
	\label{l:step1}
	It is the case that $\Pr[E_2 \cap E'_2|E_1] \geq 1 - 4 \eps_0$.
\end{Lemma}

\begin{proof}
	Assume $E_1$ has occurred and suppose $\xi_1 = \cH_\la \cap \Lambda_m$
	is known.  Define the set 
$$
\xi'_2 := \cC_{\xi_1} (
\Po_2 \cup
(\cH_\la \cap \Lambda_{n} \setminus \Lambda_{m+1} )),
$$
which can be thought of as an intermediate stage  between $\xi_1$ and $\xi_2$.

	We consider (just for the proof of this lemma)
	a second algorithm that goes in the reverse direction
(inward from the boundary of $\Lambda_n$ rather than outward
from 
	$\Lambda_m$).

Recalling that we take $n$ to be a multiple of $2m$, divide
$T(n,m)$  into cubes of side $2m$, denoted $Q_1,\ldots,Q_{t_3}$ where
	$t_3:= (n/(2m))^{d-1}$.  For $z \in T(n)$
	let $i(z)$ be the index $i \in \{1,\ldots,t_3\}$
such that
	$z +  e_1 \in Q_i$.

We create a random subset $K$ of $T(n)$ as follows.
Let $(J_z,z \in T(n))$ be Bernoulli random variables
with parameter $\eps_3$ given by \eqref{e:defeps3}, independent of each other
and of everything else.
	Then set 
	$$
	K: = \{ z \in T(n): J_z 
	{\bf 1}_{
		\{Q_{i(z)} {\rm ~is~a~seed} \}
		}
		=1
		\}.
	$$

For each 
	$y \in T(n)$, 
	define the `stick'
$
S(y):= 
	\overline{\{y,y+ e_1\}}.
$
Our choice of $\eps_3$ at \eqref{e:defeps3}
	is a lower bound for the probability
	that if we divide $S(y)$ into subcubes of side $1/(2d)$,
each subcube contains exactly  one point of $\cQ_2$,
	and the graph  $G(\cQ_2 \cap S(y))$
is connected, and given two further points in
$S(y)$, the two further points are each connected to whichever point of
$\cQ_2$ lies in the same  subcube.

Note that $K$ is constructed independently of $\xi'_2$.
Moreover, we claim that
\begin{align}
\Pr[E_2 | E_1] 
	\geq \Pr[ \overline{K} \cap \xi'_2 \neq \emptyset | E_1].
	\label{e:ForBack0}
\end{align}
	To justify this claim, define the random set
	$W:= \{i \in \{1,\ldots,t_3\}: Q_i {\rm ~ is ~ a~ seed} \}$.
	We shall condition on the sets
	$W $ and $\overline{\xi'_2} \cap T(n)$.  Given
	$\alpha \subset T(n)$ and $\beta \subset \{1,2, \ldots, t_3\}$,
	recalling the notation $\eta\|z$ 
	given shortly after \eqref{e:ndef}
	in the previous section, we compute:
	\begin{align*}
		\Pr[ E_2| 
		& \{ \overline{\xi'_2} \cap T(n) = \alpha, W = \beta\}
		\cap E_1]
		\\
		& \geq
		\Pr[ \cup_{\{z \in  \alpha :
		i(z) \in \beta \}} \{ \xi'_2 \| z
		\leftrightarrow \cH_\la \| (z + e_1) 
		{\rm ~in~} 
		G( \{\xi'_2 \| z,  \cH_\la \| (z+ e_1) \}
		\cup (\cQ_2 \cap S(z)) 
		\}
		\\
		&
		~~~~~~~~~~~~~~~
		| \{ \overline{\xi'_2} \cap T(n) = \alpha, W = \beta\}
		\cap E_1]
		\\
		& \geq
		\Pr[ \cup_{\{z \in \alpha :
		i(z) \in \beta \}} \{J_z =1\} ]
		\\
		& = \Pr[\overline{K} \cap \xi'_2 \neq \emptyset
		| \{ \overline{\xi'_2} \cap T(n)  =\alpha, W= \beta
		\} \cap E_1],
	\end{align*}
	and the claim \eqref{e:ForBack0} follows.
%

Having determined the set $K$, for the next step of the reverse algorithm
we  define the point set
$\eta_2 := \cC_{\cH_\la \cap \overline{K}} (\cH_\la \cap \Lambda_{n} 
\setminus \Lambda_{m+1})$.
Let
	$$
	U_2 :=
	( \tfrac{1}{2} + \Z)^d \cap \overline{\eta_2 \cap \Lambda_{m+2}}
	=
	\{ z \in \Delta(\Lambda_{m+1};\Lambda_{n}): \overline{\{z\}}
	\leftrightarrow \overline{K} {\rm ~in~}
	G( \cH_\la \cap  \Lambda_{n}  \setminus 
	\Lambda_{m+1}
	)
	\}.
	$$
	Also let $U_2^*$ be a maximal subset of $U_2$ having
	all vertices at $\ell^\infty$ distance at least
	5 from each other, chosen by some deterministic
	greedy algorithm  once $U_2$ is known. Then $\# U_2^*
	\geq 9^{-d} \# U_2$.

For each site $y \in U_2$, define an event  $F'_y$
as follows. 
Divide $\Lambda_{5/2}(y)$ into subcubes of side $1/(2d)$. 
Let us say that event $F'_y$
occurs if and only if
(i) there is exactly one point of $\Po_2$ in each of these
subcubes;  (ii) the graph $G(\Po_2 \cap \Lambda_{5/2}(y) )$
is connected; (iii) the first point
(in the lexicographic ordering) 
of $\xi_1$
in $\Lambda_{5/2}(y)$
is connected to the point of $\Po_2$ in its subcube, 
and (iv) the  point  $\eta_2 \| y$
is connected to the point of $\Po_2$ in its subcube.
Then with $\delta_1$ defined at \eqref{e:delta1},
conditional on $\eta_2$ and on $\xi_1$, the
events $F'_y, y\in U_2^*$ are independent
and have probability at least $\delta_1$ for each $y \in U^*_2$.
Hence by \eqref{e:tdef} we have
	\begin{align}
		\Pr[ ( \cup_{y \in U_2^*} F'_y )^c|\{\#U_2 \geq t_1\}
		\cap E_1]
		\leq (1-\delta_1)^{t_1/9^d} \leq \exp(-t_1 \delta_1/9^d) \leq
		 \eps_0.
		 \label{e:0710a}
	\end{align}

	 Next we claim that with $t_2$ given by \eqref{e:ell}
	 and $V(r,k)$ given by \eqref{e:defV}, 
	\begin{align}
		\Pr[ \{ \Lambda_m \leftrightarrow \overline{K}
		{\rm ~in~} G(\cH_\la \cap \Lambda_{n} )\}^c
		| \#V(m,n) \geq t_2] 
	\leq
	e^{-3^d \lambda t_1} \eps_0/2. 
		\label{e:LmtoK}
	\end{align}
	To see the claim, note that if $\#V(m,n) \geq t_2$
	then we can find $x_1,\ldots, x_{\lceil t_2/(2m)^{d-1}  \rceil }
	\in V(m,n)$ such that $i(x_1),\ldots, i(
	x_{\lceil t_2/(2m)^{d-1}  \rceil })$ are distinct
	and therefore the random variables 
	$J_{x_j} {\bf 1}_{\{Q_{i(x_j)} {\rm ~is~a~seed}\}}$,
	 $ 1 \leq j \leq
	\lceil t_2/(2m)^{d-1}  \rceil,$ are independent
	Bernoulli variables with parameter $\eps_2 \eps_3$, and
	hence
	\begin{align*}
		\Pr[ K \cap \{x_1,\ldots, x_{\lceil t_2/(2m)^{d-1} \rceil}\}
		= \emptyset] & \leq (1- \eps_2 \eps_3)^{\lceil
		t_2/(2m)^{d-1} \rceil}
		\\
		& \leq \exp(- \eps_2 \eps_3 t_2/(2m)^{d-1}) 
		\leq e^{-3^d \lambda t_1} \eps_0/2.  
	\end{align*}
	Also if $K \cap V(m,n) \neq \emptyset$ then $\Lambda_m 
	\leftrightarrow \overline{K}$ in $G(\cH_\la \cap \Lambda_n)$,
	and the claim follows.

By the preceding claim and the union bound,
recalling how we defined $n$ at  \eqref{e:ndef},
	we have
\begin{align}
	\label{e:Kcnct0}
	\Pr[ \{ \Lambda_m \leftrightarrow \overline{K}
	{\rm ~in~} G(\cH_\la \cap \Lambda_{n} ) \}^c] 
	\leq 
	e^{-3^d \lambda t_1}\eps_0.
\end{align}
Note that if  $\Lambda_m \leftrightarrow
\overline{K}$ in $G(\cH_\lambda \cap \Lambda_{n} )$ then also
$\Lambda_{m+1} \setminus \Lambda_{m} \leftrightarrow \overline{K}$ 
in $G(\cH_\la \cap \Lambda_{n}  )$
and therefore by \eqref{e:Kcnct0},
\begin{align*}
	\Pr[\Lambda_{m+1} \setminus \Lambda_{m} \leftrightarrow 
	\overline{K} {\rm ~in~}
	G(\cH_\la \cap \Lambda_{n}  ) ]
	\geq 1- e^{-3^d \lambda t_1}\eps_0.
\end{align*}
Hence by taking $R= \Lambda_{m+1}$ and $B= \Lambda_n$ in
Lemma \ref{l:likeGML2}, $ \Pr[\# U_2 \leq t_1] \leq \eps_0, $
and therefore using \eqref{e:0710a} and the fact that $U_2$
is independent of event $E_1$, we have
\begin{align*}
	\Pr \left[\left(\cup_{z \in U_2} F'_z\right)^c |E_1  \right]
	\leq \Pr[ \# U_2 < t_1] +
	\Pr \left[\left(\cup_{z \in U_2} F'_z\right)^c | \{ \#U_2 \geq t_1
	\} \cap
	E_1 \right]
	\leq 2 \eps_0.
\end{align*}
Suppose for some $z \in U_2$ that $F'_z \cap E_1$ occurs.
Then setting $\xi''_2 := 
\cC_{\cH_\la \cap 
\overline{K}} (\cH_\la \cap \Lambda_n \setminus \Lambda_{m+1})$,
%
we have that 
$\xi''_2 \| z$
is connected to
$\Po_2 \cap \Lambda_{1/2}(z)$ and hence
$\xi''_2 \| z
\in \xi'_2$,
so that $\xi'_2 \cap \overline{K} \neq \emptyset$. Therefore
by \eqref{e:ForBack0}, 
we have:
\begin{align*}
	\Pr[E_2|E_1] \geq \Pr[ \overline{K} \cap \xi'_2 \neq \emptyset|E_1]
	\geq \Pr[\cup_{z \in U_2} F'_z |E_1]
	\geq 1- 2 \eps_0.
\end{align*}
Similarly, $\Pr[E'_2|E_1] \geq 1 - 2 \eps_0$, and hence by the
union bound $\Pr[E_2 \cap E'_2|E_1] \geq 1-4 \eps_0$, as claimed.
\end{proof}

\section{Algorithm: Third step}
\label{s:step2}

Suppose the first two steps  were successful, i.e. $E_1 \cap E_2
\cap E'_2$ occurs.
Define the sets
\begin{align*}
	T_3(n) := \{
		x \in ( \tfrac{1}{2} +\Z)^d:
	x \cdot e_1 =n- \tfrac12, 
	0 >  x \cdot e_j  > -  n,
	{\rm ~for~} 
	j=2,3,\ldots,d
\} ;
\nonumber
	\\
	T_3(n,m)  := 
	\{	x \in \R^d : n \leq x \cdot e_1 \leq n +m, 
	0 \geq  x \cdot e_j \geq  -n, {\rm ~for~} j = 2,3, \ldots, d
	 \} .
\end{align*}
By the success of the previous steps, we have identified 
a seed $\Lambda_m(v_1)$ that is contained in $T(n,m) \cap \overline{\xi_2}$,
where $\xi_2$
is the current sub-cluster at the origin, as 
 defined at \eqref{e:defxi1} in the previous step.
 We have also identified a seed $\Lambda_m(v_2)$ 
contained in $T_2(n,m) \cap \overline{\xi_2}$.

Starting from the seed centred at $v_1$,
for the next step  of the algorithm we  see if
it is connected via the box
$\Lambda_{n}(v_1)$ 
to a further seed contained in $v_1 + T_3(n,m)$
(we call this {\em exploring to the right}).  More precisely, set
$$
\xi_3 := \cC_{\xi_2}(\Po_3 \cup \cQ_3 \cup 
(\cH_\la \cap [(\Lambda_n(v_1)  \setminus \overline{
	(\xi_2)^+_{\Lambda_n(v_1)} } \, )
 \cup  (v_1 + T_3(n,m))])) ,
$$
where $\Po_3$ and $\cQ_3$
are  homogeneous Poisson processes
of intensity $\eps_1$ in $\Lambda_{n+1}(v_1)$,
independent of $\cH_\la$, $\Po_1$, $\Po_2$, $\cQ_2$ and
each other.

Recalling the notation $\eta_B^{++}$ defined at \eqref{e:++def},
observe that 
 \begin{align}
 \overline{(\xi_2)^{++}_{\Lambda_n(v_1)}}  \cap ( v_1  + T_3(n))
 = \emptyset,
	 \label{e:Tnew}
 \end{align}
 since $x \cdot e_1 \leq v_1 \cdot e_1 + m$ for all
 $x \in \xi_2$, whereas $y \cdot e_1 = v_1 \cdot e_1 + n - \frac12$
 for $y \in v_1 + T_3(n)$.

We denote by $E_3$ the event that 
there is a seed contained in 
$\overline{\xi_3} \cap (v_1 + T_3(n,m))$.
If $E_3$ occurs, then we say this step of the algorithm  is successful, and
denote the centre of the first such seed by $v_3$.
\begin{Lemma}
	\label{l:step2}
	It is the  case that $\Pr[E_3|E_1 \cap E_2 \cap E'_2]
	\geq 1- 2 \eps_0$.
\end{Lemma}

\begin{proof}
	The proof is similar to that of Lemma \ref{l:step1}.
We divide the slab $v_1 + T_3(n,m)$ into
cubes of side $2m$, denoted $\tQ_1,\ldots,\tQ_{t_3}$, where
	$t_3:= (n/(2m))^{d-1}$ as before.  For $z \in v_3 +  T_3(n)$
let $\tii(z)$ be the index $i$ such that $z + e_1 \in \tQ_i$.
	Again,
we consider another algorithm going in the reverse direction.
Let $(J_z,z \in T_2(n))$ be Bernoulli random variables
with parameter $\eps_3$ (defined earlier at \eqref{e:defeps3}),
independent of each other
and of everything else.
	Define the random set $$
	K_3: =
	\{ z \in v_1 + T_3(n):
	J_z 
	{\bf 1}_{\{\tQ_{\tii(z)} {\rm ~is~a~seed}\}}
	=1\}.
	$$

	Let $\xi'_3 := \cC_{\xi_2}(\Po_3 \cup [\cH_\la
	\cap (\Lambda_n(v_1)
	\setminus \overline{(\xi_2)^+_{\Lambda_n(v_1)}})])$.
	Then $K_3$ is independent of $\xi'_3$.
	As at \eqref{e:ForBack0},
	we can show by conditioning on
	$K_3$ and on $\overline{\xi'_3} \cap (v_1 +T_3(n))$ that
\begin{align}
\Pr[E_3] 
	\geq 
	\Pr[ \overline{K_3} \cap \xi'_3 \neq \emptyset],
	\label{e:ForBack}
\end{align}
where we here write just $\Pr$ for conditional probability
	given $E_1 \cap E_2 \cap E'_2$.

Having determined the set $K_3$, we set
$\eta_3 := \cC_{\cH_\la \cap \overline{K_3}}
	(\cH_\la \cap \Lambda_n(v_1) 
	\setminus \overline{(\xi_2)^+_{\Lambda_n(v_1)}})$.
Let
	\begin{align*}
		U_3 & :=
	(\tfrac12 + \Z)^d \cap \overline{\eta_3 \cap
		\overline{\Delta^2(\xi_2;\Lambda_n(v_1))}}
		\\
		& =
	\{ z \in \Delta^2(\xi_2;\Lambda_n(v_1)): \overline{\{z\}}
	\leftrightarrow \overline{K_3} {\rm ~in~}
		G( \cH_\la \cap  \Lambda_n(v_1)  \setminus 
		\overline{(\xi_2)^{+}_{\Lambda_n(v_1)}})
	\},
	\end{align*}
	and let $U_3^*$ be a subset of $U_3$ with all
	elements distant at least $5$ from each other in
	the $\ell^\infty $ norm and with $\# U_3^* \geq 9^{-d} \# U_3$.

For each site $z \in U_3$, we define an event $F''_z$ as follows.
Divide $\Lambda_{5/2}(z) \cap \Lambda_n(v_1)$ into subcubes of side $1/(2d)$. 
Let us say $F''_z$ occurs if and only if
(i) there is exactly one point of $\Po_3$ in each of these
subcubes;  (ii) the graph $G(\Po_3 \cap \Lambda_{5/2}(z) )$
is connected; (iii) the first point
(in the lexicographic ordering) 
of $\xi_2$ in $\Lambda_{5/2}(z)$
is connected to the point of $\Po_3$ in its subcube; and
(iv) the 
point  
$\eta_3\|z$
is connected to the point of
$\Po_3$ in its subcube.
Then with $\delta_1$ defined at \eqref{e:delta1},
conditional on $\eta_3$ and on $\xi_2$, the events
$F''_z$, $z \in U_3^*$ are independent and have probability
at least $\delta_1$ for each $z \in U_3^*$. Hence using 
\eqref{e:tdef}, we have
	$\Pr[ (\cup_{z \in U_3^*} F''_z)^c|\#U_3
	\geq t_1] \leq  \eps_0$ as in
	\eqref{e:0710a}.

	Now
	set
	\begin{align*}
		V'(m,n) := \{x \in v_1 + T_3(n): \Lambda_m(v_1) \leftrightarrow
		\Lambda_{1/2}(x) {\rm ~in~} G(\cH_\la \cap 
		\Lambda_{n}(v_1) )
		\}.
	\end{align*}
	Let $t_2$ be given by \eqref{e:ell}.
	Then as at \eqref{e:LmtoK},
	\begin{align*}
		\Pr[ \{ \Lambda_m(v_1) \leftrightarrow \overline{K_3}
		{\rm ~in~} G(\cH_\la \cap \Lambda_{n}(v_1) )\}^c
		| \#V'(m,n) \geq t_2] 
	\leq
	 e^{-3^d \lambda t_1} \eps_0/2. 
	\end{align*}
By \eqref{e:ndef}
and symmetry,
	\begin{align*}
	\Pr[\# V'(m,n) < t_2 ] 
	= \Pr[\# V(m,n) < t_2 ] 
		\leq
		e^{-3^d \lambda t_1} \eps_0/2, 
	\end{align*}
	and therefore,  since $\Lambda_m(v_1) \subset \overline{\xi_2}$
	and \eqref{e:Tnew} holds,
\begin{align*}
	\Pr[	\Delta^2 ( \xi_2;\Lambda_n(v_1)) \leftrightarrow 
	\overline{K_3} {\rm ~in~}
	G(\cH_\la \cap \Lambda_{n}(v_1)  \setminus
	\overline{(\xi_2)_{\Lambda_n(v_1)}^+}) ]
	\\
	\geq
	\Pr[ \Lambda_m(v_1) \leftrightarrow \overline{K_3}
	{\rm ~in~} G(\cH_\la \cap \Lambda_{n}(v_1) )] 
	\geq 1- 
		e^{-3^d \lambda t_1} \eps_0. 
\end{align*}
Hence, by taking $R  = (\xi'_3)^+_{\Lambda_n(v_1)}$ and
$B= \Lambda_n(v_1)$ in Lemma \ref{l:likeGML2}
(which is applicable by \eqref{e:Tnew}), 
$
\Pr[\# U_3 \leq t_1] \leq \eps_0,
$
and therefore
\begin{align*}
	\Pr \left[\left(\cup_{z \in U_3} F''_z\right)^c \right]
	\leq \Pr[ \# U_3 \leq t_1] +
	\Pr \left[\left(\cup_{z \in U_3} F''_z\right)^c | \#U_3 \geq t_1 \right]
	\leq 2 \eps_0.
\end{align*}
Suppose for some $z \in U_3$ that $F''_z$ occurs.  Then setting
$\xi''_3 := \cC_{\cH_\la \cap \overline{K_3} }
(\cH_\la \cap \Lambda_n \setminus \overline{(\xi_2)^{+}_{\Lambda_n}})$
we have that $\xi''_3 \| z$ is connected to
$\Po_3 \cap \Lambda_{1/2}(z)$ and hence $\xi''_3 \| z \in \xi'_3$,
so that $\xi'_3 \cap \overline{K_3} \neq \emptyset$. Therefore
by \eqref{e:ForBack} we have:
\begin{align*}
	\Pr[E_3] \geq \Pr[ \overline{K_3} \cap \xi'_3 \neq \emptyset]
	\geq \Pr[\cup_{z \in U_3} F''_z]
	\geq 1- 2 \eps_0,
\end{align*}
as required.
\end{proof}

\section{Continuing the algorithm}
\label{s:continuing}

The rest of the argument proceeds along similar lines
to \cite{GM90}; we give enough details to show that
our earlier choices of $\eps_0$ and other constants in
Section \ref{s:step0} work.

We are eventually going to show that our process dominates
a supercritical oriented site percolation process on $\Z_+^2$,
where $\Z^2_+ := (\N \cup \{0\})^2$.
For each $(i,j) \in \Z_+^2$, we
shall determine the status of site $(i,j)$ (open/closed)
on the basis of the previous history of our algorithm,
along with examination of $\cH_\la$ and our sprinkling
Poisson processes in certain new parts of
space in the vicinity of $4N (ie_1 + je_2)$, where we set $N :=n+m$. 
Set $B_{(i,j)} := \Lambda_{N}(4N (ie_1 + j e_2))$. We refer to all such
$B_{(i,j)}$ as a {\em site-boxes}. We also refer to
$\Lambda_N(4N(ie_1+je_2) + 2Ne_1)$ and $\Lambda_N(4N(ie_1 + je_2) + 2Nn e_2)$
as {\em bond-boxes}.

\begin{figure}[!h]
\center
\includegraphics[width=10cm]{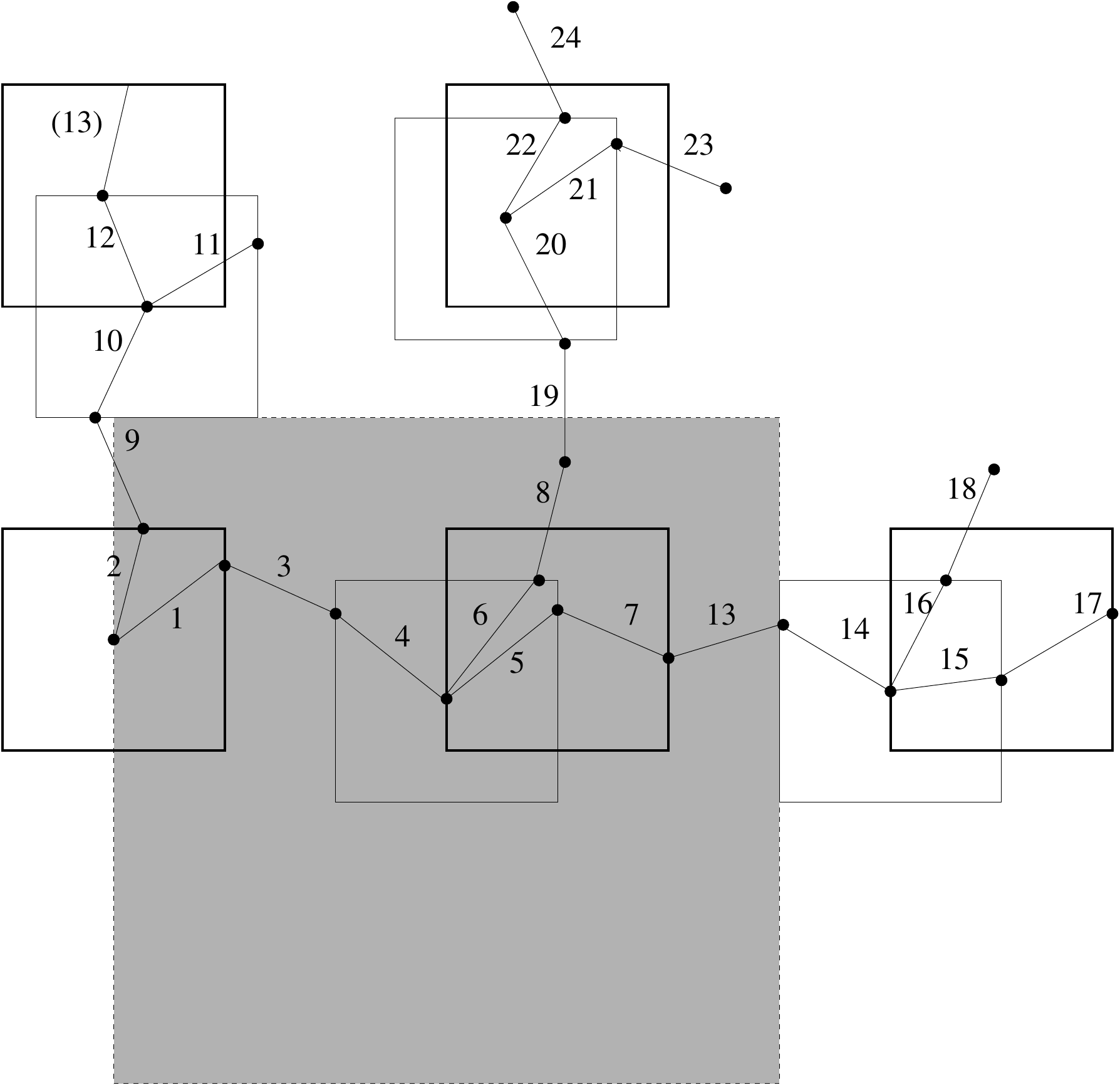}
	 \caption{\label{fig2}The bold squares represent
	  site-boxes. The dots represent seeds,
	 with
	  seed $v_i$ at the end (further from the root)
	  of $i$th edge. The step starting from $v_{12}$ is assumed
	  to be unsuccessful so
	  $(0,1)$ is closed;
	  all of the other  steps shown are successful. All Poisson processes involved in stage $(1,0)$ are inside the shaded region.}
\end{figure}

Some of the site-boxes are shown in Figure \ref{fig2}.
Each bond-box lies between two neighbouring site-boxes.

The algorithm proceeds in {\em stages}, each
stage being associated with a vertex of $\Z_+^2$. We consider
these vertices in order of increasing $\ell^1$-distance
from the origin; among vertices the same $\ell^1$ distance
from the origin we consider them in decreasing lexicographical
order (so our ordering starts with $(0,0), (1,0), (0,1),
(2,0), (1,1), (0,2), (3,0), \ldots$).
Each stage is broken up into at most five {\em steps}.
Given $n \in \N$,
in the $(n+1)$-st step we start with a set $\xi_{n}$
arising at the end of the previous step, and
update this to a new set $\xi_{n+1}$.

The first two steps which we already discussed in Sections \ref{s:step0} and
\ref{s:step1}, 
taken together, will be considered
as a whole to be the first stage of the algorithm,
associated with the vertex $(0,0)$ of $\Z^2_+$.
If both of these  steps  were successful, i.e. if events $E_1, E_2 $ and 
$E'_2$ all occur, then
we deem the first stage as a whole (associated with $(0,0)$ in
$\Z_+^2$) to be successful, and site $(0,0)$ to be {\em open} (otherwise
$(0,0)$ is closed).

In general, we say the stage associated with site $(i,j) \in \Z_+^2$
is successful (and $(i,j)$ is open) if
we succeed in finding seeds centred in the bond-boxes just above and
just to the right of $B_{(i,j)}$, included in the current value
of the sub-cluster at the origin after the completion  of this stage.
For the $(0,0)$ stage the seeds created are centred right on the boundary 
of their bond-boxes, but we take all boxes to be closed so this
counts.

Suppose the first stage of the algorithm (associated
with vertex $(0,0)$ of $\Z^+_2$) was successful
(so the first and second steps were both successful,
i.e. $E_1 \cap E_2 \cap E'_2$ occurs).
Then the third step is as already described in Section
\ref{s:step2}, and is deemed successful 
if $E_3$ occurs.

If $E_3$ occurs, then there is a 
 further seed contained in
$\overline{\xi_3} \cap
 (v_1 + T_3(n,m))$.
We denote the centre of the first such seed in the lexicographic ordering
by $v_3$.

Repeating this, if the previous step is successful
we continue to explore to the right
with a further similar step starting from the new seed centred at $v_3$
and extending $\xi_3$ to a larger sub-cluster $\xi_4$.
In this fourth step, we steer all other coordinates
other than the first one towards zero (we shall
give more detail on the steering in the next step).
If this step is also successful we arrive at
a seed contained in  $\overline{\xi_4}$ and with its 
centre (denoted $v_4$) lying  in the right
face of $\Lambda_N(v_3)$ and in the left face 
of $B_{(1,0)}$ (see Figure \ref{fig2}).

Suppose the last step described is also successful.
We need to `branch out' from the seed $\Lambda_m(v_4)$,
seeking further connections from $\xi_4$ to seeds adjacent
to both  the
easterly and the northerly face of $\Lambda_n(v_4)$,
and choosing which orthant of the face by
an appropriate steering mechanism in both cases.
For this, we follow a similar procedure
to the second step as described in Section \ref{s:step1},
but now starting from the seed $\Lambda_m(v_4)$.
To go into more detail, 
for $s \in \R\setminus \{0\}$ set $\sgn(s) = s/|s|$, and set 
$\sgn(0) := 1$.
Define
\begin{align*}
	T_5(n) := \{
		x \in ( \tfrac{1}{2} +\Z)^d:
	x \cdot e_1 =n- \tfrac12, 0 > (\sgn (v_4 \cdot e_j)) x \cdot e_j > - n
	{\rm ~for~} 
	2 \leq j \leq d
\} ;
\nonumber
	\\
	T_5(n,m)  := 
	\{	x \in \R^d :
	n \leq x \cdot e_1 \leq n +m, 
	0 \geq (\sgn (v_4 \cdot e_j)) x \cdot e_j \geq - n
	{\rm ~for~} 
	2 \leq j \leq d
	 \} ,
\end{align*}
and
\begin{align*}
	T_6(n) := \{ &
		x \in ( \tfrac{1}{2} +\Z)^d:
	x \cdot e_2 =n- \tfrac12, 
	0 <  x \cdot e_1 <  n
	\\
	&  0 > (\sgn (v_4 \cdot e_j)) x \cdot e_j > - n
	{\rm ~for~} 
	3 \leq j \leq d
\} ;
\nonumber
	\\
	T_6(n,m)  := 
	\{	& x \in \R^d :
	n \leq x \cdot e_2 \leq n +m, 
	0  \leq x \cdot e_1 \leq n, 
	\\
	& 0 \geq (\sgn (v_4 \cdot e_j)) x \cdot e_j \geq - n
	{\rm ~for~} 
	3 \leq j \leq d
	 \} .
\end{align*}
We then set 
$$
\xi_5 := \cC_{\xi_4}(\Po_5 \cup \cQ_5 \cup
(\cH_\la \cap [(\Lambda_{n}(v_4) \setminus
\overline{(\xi_4)^+_{\Lambda_{n}(v_4) }} \: )
\cup (v_4 + T_5(n,m)) \cup (v_4 + T_6(n,m))])) ,
$$
where $\Po_5$ and $\cQ_5$ are independent 
homogeneous Poisson processes
of intensity $\eps_1$ in $\Lambda_{n+1}(v_4)$,
independent of all Poisson processes mentioned already.

Let $E_5 $ be the event that there is  further
seed contained in
$\overline{\xi_5} \cap (v_4 + T_5(n,m))$; if
this occurs,
denote the centre of the first such seed by $v_5$.
Let $E_6 $ be the event that there is  further
seed contained in
$\overline{\xi_5} \cap (v_4 + T_6(n,m))$; if
this occurs,
 denote the centre of the first such seed by $v_6$.

\begin{figure}[!h]
\center
\includegraphics[width=10cm]{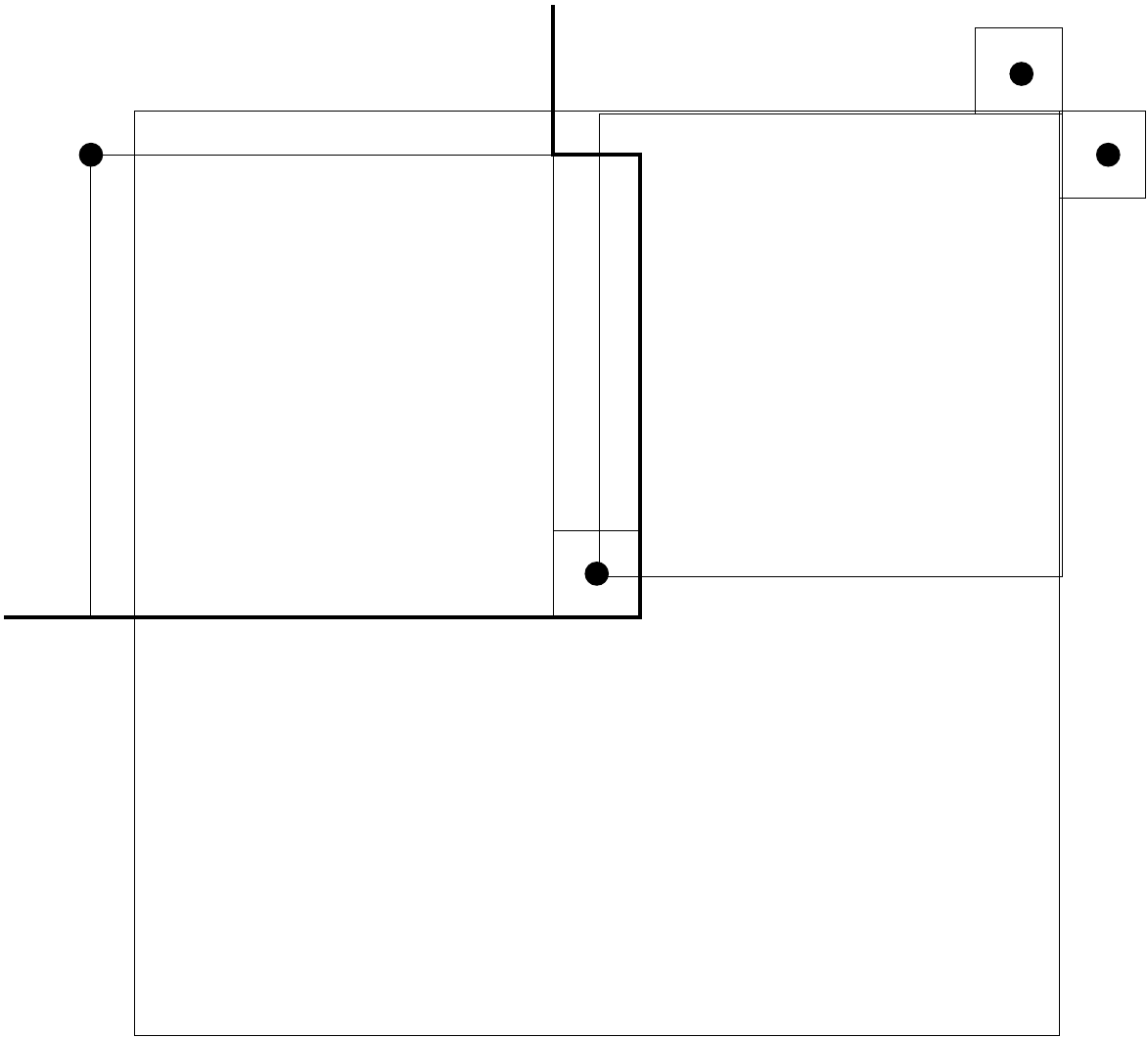}
	 \caption{\label{fig1}The largest square represents
	 $\Lambda_{n}(v_4)$; the medium-sized squares are
	 parts of $\Lambda_n(v_3)$ and of $\Lambda_n(v_4)$,
	 and the small squares (from left to right) are $\Lambda_m(v_3)$,
	 $\Lambda_m(v_4)$, $\Lambda_m(v_6)$ and $\Lambda_m(v_5)$.
	 The thick polygonal line represents part of the boundary
	 of the region possibly explored in creating $\xi_4$, illustrating
	 \eqref{e:Tnew2}.}
\end{figure}

 Note that in our definition of both $T_5(n,m)$ and $T_6(n,m)$
 we steer 
 the 3rd, 4th, $\ldots$ $d$th coordinates towards zero.  
 With $T_5(n,m)$ we likewise steer the 2nd  coordinate towards zero.
 With $T_6(n,m)$,  
we steer the first coordinate to the right.
This, and the assumption $m \geq 9$,
ensures that  the new seed we are trying to create next to the
north face of $\Lambda_n(v_4)$ is adjacent to a previously
unexplored region, i.e.
\begin{align}
 \overline{(\xi_4)^{++}_{\Lambda_n(v_4)}}  \cap ( v_4  + T_6(n))
 = \emptyset,
	 \label{e:Tnew2}
\end{align}
as illustrated in Figure  \ref{fig1}.
Since $v_4$ lies in the left half of $B_{(1,0)}$, this steering also
ensures that
 $v_6$ will lie within
$\ell^\infty$-distance $N$ of the line $\{4Ne_1 + u e_2: u \in \R\}$.


 If all the earlier steps were successful, then similarly
 to Lemma \ref{l:step2} both $E_5$ and $E_6$ have probability
 at least $1 - 2 \eps_0$, so by the union bound $E_5 \cap E_6$
 occurs with probability at least $1 - 4 \eps_0$. If
 $E_5 \cap E_6$ occurs we say this step is successful. Thus
 in each of the successful `branching' steps, two new seeds are created.

Assuming both of these steps were successful, we proceed
with a further step to the right from the new seed we found
next to the right face of $\Lambda_n(v_4)$; if successful,
this will give
us a new seed with centre in the bond-box $\Lambda_{N}(6N e_1)$
(in fact, on its boundary). 

The centre of the new seed we found next to the upper boundary 
of $\Lambda_n(v_4)$
(denoted $v_6$)
might itself lie in the bond-box
$\Lambda_N(4Ne_1 + 2N e_2)$; if not
(as in the case illustrated in Figure \ref{fig2}),
we proceed with a further step upward in which we seek a further new seed
adjacent to the upper face of $\Lambda_n(v_6)$, and if this is
successful, then  the centre of the further new seed will lie in the bond-box
$\Lambda_N(4Ne_1 + 2N e_2)$.

If all of the steps just described are successful, we declare 
this stage
to be  successful 
and deem $(1,0)$ to be open.  What we have achieved in this stage,
starting from a seed centred in the bond-box to the left of the site-box
$B_{(0,1)}$, is to extend our sub-cluster to include a seed in 
the bond-box centred just  to the right of $B_{(1,0)}$, and
 a further seed centred in the bond-box  just above $B_{(1,0)}$.

After this stage, we go on to
$(0,1)$ for the next stage.

Let us describe how our algorithm proceeds at a generic stage,
associated with vertex $(i,j) \in \Z_{+}^2$.
First suppose there is an oriented path of open sites
in $\Z_+^2$ from $(0,0)$ to $(i-1,j)$. Suppose also that
$k$ steps have been completed so far, and $\ell +1$ new seeds
(with centres denoted $v_0, v_1,\ldots,v_\ell$, where $v_0 := \0$) have
been identified so far. Here $\ell \leq 2k$ since each step involves
identifying at most two new seeds; for example, the seeds $v_1$ and $v_2$
that were identified at the end of the second step of the algorithm.

Then the current  sub-cluster at the origin $\xi_k$ is such that
$\overline{\xi_k}$  contains a seed with its centre, denoted $v_{\ell'}$,
say, where $\ell' \leq \ell$,
in the left half of the bond-box to the left of $B_{(i,j)}$. 
We then take two steps to the right, to
get a seed with its centre in the left half of the box
$B_{(i,j)}$. Let us denote the centre of this new seed  by $v_{\ell+2}$,
and new extended sub-cluster at the origin by $\xi_{k+2}$.
By following analogous steps to those described
in Section \ref{s:step1},
but
starting from the seed $\Lambda_m(v_{\ell +2})$, we extend
$\xi_{k+2}$ to a further sub-cluster containing  seeds with centres
adjacent to both 
 the right face and the upper face of $\Lambda_n(v_{k+2})$,
 denoted $v_{\ell +3}$ and $v_{\ell +4}$.

 Taking one further step to the right from $v_{\ell +3}$
 we arrive at a seed in the bond-box immediately to the right
 of $B_{(i,j)}$, and taking a northward step 
 (if needed) from $v_{k+4}$ we arrive at a seed in the
 bond-box just above $B_{(i,j)}$. If all of these steps
 are  successful we declare $(i,j)$ to be open; otherwise
 we declare it to be closed.

Now suppose there is
no  oriented path of open sites
in $\Z_+^2$ from $(0,0)$ to $(i-1,j)$, but there is
an  oriented path of open sites
in $\Z_+^2$ from $(0,0)$ to $(i,j-1)$.
Then the current sub-cluster at the origin (again denoted $\xi_k$) includes
a seed with its centre (again denoted $v_{\ell'}$)
in the bond-box just below
$B_{(i,j)}$. We then take two upward steps
from $v_{\ell'}$
to get a seed with its centre
in the bottom half of the box
$B_{(i,j)}$. Let us denote the centre of this seed  by $v_{\ell+2}$ and the
new extended sub-cluster at the origin by $\xi_{k +2}$.
By following the step in Section \ref{s:step1}
but starting from the seed $\Lambda_m(v_{\ell +2})$, we extend
$\xi_{k+2}$ to a further sub-cluster containing  seeds adjacent to both 
 the right face and the upper face of $\Lambda_n(v_{\ell +2})$,
 with centres
 denoted $v_{\ell +3}$ and $v_{\ell +4}$ respectively.

  Taking a rightward step 
 (if needed) from $v_{\ell +3}$ we arrive at a seed in the
 bond-box just to the right of  $B_{(i,j)}$,  
 and taking one further upward step from $v_{\ell +4}$
 we arrive at a seed in the bond-box immediately above
  $B_{(i,j)}$.
 If all of these steps
 are  successful we declare $(i,j)$ to be open (otherwise
 we say it is closed).

 If there is no
  oriented path of open sites
in $\Z_+^2$ from $(0,0)$ either to $(i,j-1)$ or to $(i-1,j)$
(in which case we say $(i,j)$ is {\em inaccessible}),
then we declare $(i,j)$ to be open with probability $1- 20 \eps_0$
and otherwise closed, without extending the set $\xi_k$.
As soon as we have determined the open/closed status of $(i,j)$
we go on to the next site of $\Z_+^2$.

Thus if $(i,j)$ is open, and there is an oriented path in
$\Z_+^2$ of open
sites from $(0,0)$ to $(i,j)$, then the current
sub-cluster at the origin (after the stage associated with
$(i,j)$ is completed), denoted $\xi_m$ say,
is such that $\overline{\xi_m}$
contains both  a seed centred in the bond-box
just above $B_{(i,j)}$, and a seed centred in the bond-box
just to the right of $B_{(i,j)}$.

\begin{proof}[Proof of Theorem \ref{t:GMphi}]
	As discussed at the start of Section \ref{s:step0},
	it suffices to prove that there exists some $M$ such
	that $\mu \geq \lambda_c(\phi,S_M)$.
	We shall do so for $M = 4N = 4(n+m)$.

	In our algorithm, each stage takes at most five steps, and each step
	involves
sprinkling at most two Poisson point processes of intensity $\eps_1$ in 
	a box $\Lambda_{n+1}(v_\ell)$, where $v_\ell$ is the
	centre of the starting box for that step;
	for example the Poisson process $\Po_1$ 
	in the very first step described in Section
	\ref{s:step0} (in fact this step involves sprinkling on the smaller
	box $\Lambda_{1/2}$ rather than $\Lambda_{n+1}$),
	the Poisson processes $\Po_2$ and $\cQ_2$ 
	in the second  step from Section \ref{s:step1}, and
the  processes $\Po_3$ and $\cQ_3$  in the  the third step
from Section \ref{s:step2}.

To ensure independence of the different Poisson processes we use,
we shall use different sprinkling Poisson processes for
each of the (at most) five steps within a particular stage,
and also for stages associated with distinct but adjacent
	(in the $\ell^\infty$ sense)
	sites of $\Z_+^2$. 

	For $(i,j) \in \Z_+^2$ and $k =1,2,3,4,5$ let
	$\Po_{i,j,k}$ and $\cQ_{i,j,k}$
	be the Poisson process of sprinkled points that we use
	for the $k$th step of stage $(i,j)$ of the algorithm; for example the 
	$\Po_3$ and $\cQ_3$ that we used in the 
	first step of stage $(1,0)$ in Section
	\ref{s:step2} become $\Po_{1,0,1}$ and $\cQ_{1,0,1}$
	in this notation. 

	All Poisson processes involved in stage $(1,0)$ are
	inside the region $\Lambda_{3N}(3Ne_1- e_2) $ (the shaded
	region in Figure \ref{fig2}). More generally,
	for any $(i,j,k)$ the Poisson processes $\Po_{i,j,k}$
	and $\cQ_{i,j,k}$ are inside the region $\Lambda_{3N}((4i-1)Ne_1
	+ (4j-1)N e_2)$.
	Thus, for any pair of sites
	(say $(i,j)$ and $(i',j')$)
	of $\Z_+^2$ at $\ell^\infty$ 
	distance $2$ or more from each other, the
	sprinkled Poisson processes involved in the steps that
	are  part of stage $(i,j)$ are in a disjoint region from 
	those involved in steps that are part of stage $(i',j')$.

	Let $k \in \{1,\ldots,5\} $.
	The Poisson processes $\Po_{1,0,k},$ $ \Po_{3,0,k},$ $ \Po_{1,2,k},$ $ \Po_{5,0,k},$ $ \Po_{3,2,k},$ $
	\Po_{1,4,k}, \ldots$
	are on disjoint regions of $\R^d$ and can be considered as part of a single
	homogeneous Poisson process of intensity $\eps_1$ on $\R^d$.
	Likewise
	 $\cQ_{1,0,k},$ $ \cQ_{3,0,k},$ $ \cQ_{1,2,k},$ $ \cQ_{5,0,k},$ $ \cQ_{3,2,k},$ $
	\cQ_{1,4,k}, \ldots$ can be considered as part of a single
	homogeneous Poisson process of intensity $\eps_1$ on $\R^d$.
	Similarly we
	can take the Poisson processes $\Po_{i,j,k}$ for all $(i,j)$ having the
	same parity as $(0,0)$ to be all part of a single homogeneous Poisson process
	of intensity $\eps_1$, and likewise for all $(i,j)$ having the 
	same parity as $(0,1)$, and for all $(i,j)$ having the same parity as $(1,1)$ (and similarly with the $\cQ_{i,j,k}$'s).

Therefore all of the Poisson processes of sprinkled points in the algorithm
	can be considered to be part of a total of 
	$2 \times 5 \times 4 = 40$ independent
	homogeneous Poisson processes 
	of intensity $\eps_1$ in $\R^d$.
Since we took $\eps_1 = (\mu-\la)/40$, 
this means the union of 
	all of these Poisson processes
is stochastically dominated by a homogeneous Poisson process
	of intensity $\mu - \lambda$, independent
	of $\cH_\la$, and thus at all stages the
	set $\xi_n$ can be considered to be a subset of $\cH_{\mu}
	\cap S_{4N}$, since the steering means that for
	each  $\ell \in \N$ and $j=3,4,\ldots,d$
	the $j$-coordinate of  $v_\ell$
	stays in the range $[-n,n]$, and the
	 box $\Lambda_{n+m}(v_\ell)$ lies within
	$S_{4N}$.

	By Lemmas \ref{l:step1}
	and \ref{l:step2}, and similar arguments subsequently,
	each step of the algorithm is successful with probability
	at least $1- 4 \eps_0$. There are at most five steps in
	each stage. Hence by the union bound, if
	$(i,j) \in \Z_+^2$ is such that there is an oriented path of
	open sites from $(0,0)$ either to $(i-1,j)$ or to $(i,j-1)$,
	then the stage
	associated with $(i,j)$  is successful in all its steps
	with probability at least $1- 20 \eps_0$.
	Since we also take $(i,j)$ to be open with probability $1-20 \eps_0$
	if it is inaccessible, the conditional
	distribution of the random field $({\bf 1}_{\{(i,j) {\rm ~is~
	open}\}}, (i,j) \in \Z_+^2)$, given $E_1$ occurs,
	stochastically dominates
	an independent Bernoulli random field with parameter $1- 20 \eps_0$.
	Since we took  $\eps_0 = 1/9999$, we have $1-20 \eps_0 > 80/81$.
	The critical probability $p_c$
for oriented site percolation on 
$\Z_+^2$ 
satisfies
$p_c \leq 80/81$ 
\cite[p. 1026]{Dur84}, and therefore the probability
that there is an infinite oriented path of open sites starting at $(0,0)$ is 
strictly positive.

If there is an infinite oriented path in $\Z_+^2$  of open sites starting from
$(0,0)$, then
	the union $\cup_{i=1}^\infty \xi_i$ is
	an infinite set of vertices of $\cH^\0_{\mu} \cap S_{4N}$ that includes
	the origin and induces a connected subgraph of $G(\cH^{\0}_\mu)$,
	so that  
		%
		$\0 \leftrightarrow \infty {\rm ~in~} G((\cH_\mu \cap S_{4N }) \cup \{\0\} )
		$.
	Therefore  $\mu \geq \lambda_c(S_{4N})$ as required.
%
\end{proof}

\section{The giant component}
\label{s:Giant}

For any finite graph $G$, let $L_j(G)$  denote the order of
its $j$th-largest component, that is, the $j$th-largest of the 
orders of its components, or zero if it has fewer than $j$
components.
Theorem \ref{t:GMphi}
is used in the proof of the following result on the {\em giant component}
phenomenon for the random connection model in a finite box (which is also known
as the {\em soft random geometric graph}).

%
%

\begin{Theorem} \cite{KP26}
	\label{thm1}
Suppose $d \geq 3$.
	Suppose $\sup\{r >0: \phi(r) >0\} \in (0,\infty)$
	and $\la > \la_c(\phi)$. Then
	as $s \to \infty$ we have that
	\begin{align}
		(2s)^{-d}L_{1}(G(\cH_{\la,s},\phi)) \toL \la \theta(\phi,\la,\R^d).
	\end{align}
\end{Theorem}


Let $\cC_\infty(\cH_\la)$ be the set of vertices $x \in \cH_\la$ such that
$x \leftrightarrow \infty$ in $G(\cH_\la,\phi)$ (i.e., the infinite cluster
of this graph), and let $R_\la
:= \sup\{r \geq 0: \cC_\infty(\cH_\la) \cap \Lambda_r = \emptyset \}$.
K\"upper  \cite[Theorem 5.1]{K25} proves a
weaker version of Theorem \ref{thm1}, 
namely 
\begin{align}
	\lim_{s \to \infty} \Big( \frac{\E[L_1(G(\cH_{\lambda,s},\phi))]}{
		(2s)^d} \Big) = \lambda  \theta(\phi,\lambda,\R^d),
		\label{e:nicgiant}
\end{align}
under the assumption (not proved in \cite{K25}) that
$\Pr[ R_\la > r] = O(r^{-\delta})$
as $r \to \infty$ for some $\delta >0$. 
Using our Theorem \ref{t:GMphi}, we
shall show in Corollary \ref{c:tail} below that
this assumption is justified, and in fact for all
$\la > \la_c(\phi)$ we have
a stronger tail estimate for $R_\la$, namely 
$\Pr[R_\la >r] \leq \exp(-c r^{d-2})$ for some constant $c >0$.

In \cite{KP26}, the proof of
Theorem \ref{thm1} is completed, using K\"upper's result
\eqref{e:nicgiant} as a starting point. 
The case $d=2$ was considered earlier in \cite{P22}.

\begin{Corollary}
	\label{c:tail}
	Suppose $d \geq 3$ and $\lambda > \lambda_c(\phi,\R^d)$.
	Then
	\begin{align}
		\limsup_{r \to \infty} r^{2-d}
		\left(	\log \Pr[\{\Lambda_r \leftrightarrow \infty
		~{\rm in}~ G(\cH_\la,\phi) \}^c ] \right) < 0.
		\label{e:expdecay}
	\end{align}
\end{Corollary}
\begin{proof}
	Recall that $\rho_\phi$ denotes the range of $\phi$.
	Using Theorem \ref{t:GMphi}, choose $M \in (2 \rho_\phi,\infty)$
	such that $\lambda > \lambda_c(\phi,S_M)$. 
	Then given $r >M$, we can find $x_1,x_2,\dots x_{\lfloor r/M
	\rfloor^{d-2}}$ such that $x_i \cdot e_1 = x_i \cdot e_2 = 0$
	for all $i \in \{1,\ldots, \lfloor r/M \rfloor^{d-2}\}$
	 and $\Lambda_{M}(x_1) , \ldots, \Lambda_{M}(x_{\lfloor r/M \rfloor^{d-2}}
	 ) $ are subsets of $\Lambda_r$ with disjoint interiors.
	Then $x_1 + S_{2M}, \ldots, x_{\lfloor r/M \rfloor^{d-2}} + S_{2M}$
	have disjoint interiors, and hence
	\begin{align*} 
		\Pr[\{\Lambda_r \leftrightarrow
		\infty {\rm ~ in~}
		G(\cH_\la)\}^c]
		& \leq \Pr[ \cap_{i=1}^{\lfloor r/M\rfloor^{d-2}}  
		\{ 
			 \Lambda_{\rho_\phi}(x_i) \leftrightarrow  \infty
		{\rm ~in~} G(\cH_\la \cap (x_i + S_{2M} ) ) \}^c]
		\\
		& \leq 
		\prod_{i=1}^{\lfloor r/M\rfloor^{d-2}}  
		(1- 
		\Pr[ \0 \leftrightarrow \infty ~{\rm in}~ G(\cH_\la \cap S_M)]).
	\end{align*}
	Since $\theta(\phi,\lambda,S_M) >0$,
	\eqref{e:expdecay}
	follows.
\end{proof}

{\bf Acknowledgement.} I thank Ivailo Hartarsky for drawing 
my attention to \cite{FH25}, and the referee for some helpful comments.

\end{document}